\documentclass[10pt]{article} 
\usepackage[margin=1in]{geometry}

 \usepackage{graphicx}
\usepackage{listings}
\usepackage{geometry}

\usepackage{amsmath,amsthm,amsfonts,amssymb,amscd}
\usepackage[colorlinks=true, pdfstartview=FitV, linkcolor=blue,citecolor=blue, urlcolor=blue]{hyperref}
\usepackage[usenames,dvipsnames]{color}
\usepackage{times}
\usepackage{comment}

\usepackage{mathtools}
\usepackage{mathabx}
\usepackage{float}
\usepackage{makeidx}
\usepackage{stmaryrd}
\usepackage{mathrsfs}
\usepackage{emptypage}
\usepackage{xfrac}
\usepackage{bbm}

\renewcommand{\epsilon}{\varepsilon}

\usepackage{graphicx}
\newcommand{\p}{\ensuremath{\partial}}

\definecolor{labelkey}{rgb}{0,0,1}

\def\les{\lesssim}

\def\eps{\varepsilon}
\renewcommand*{\div}{\ensuremath{\mathrm{div\,}}}

\renewcommand*{\tilde}{\widetilde}
\renewcommand*{\hat}{\widehat}
\renewcommand*{\bar}{\overline}

\newtheorem{theorem}{Theorem}[section]

\theoremstyle{definition}

\newtheorem{remark}[theorem]{Remark}

\numberwithin{equation}{section}

\def\p{\partial}

\def\f1r{{\frac{1}{r}}  }
\def\p{\partial}

\def\f1r{{\frac{1}{r}}  }

\allowdisplaybreaks[4]

\title{A note on the existence of self-similar profiles of  the hydrodynamic formulation of
the focusing nonlinear Schr\"odinger equation}
  \author{G. Cao-Labora, J. G\'omez-Serrano, J. Shi, G. Staffilani}

\date{} %
\begin{document}

\maketitle
 \begin{abstract}
 After performing the Madelung transformation, the nonlinear Schr\"odinger equation is transformed into a hydrodynamic equation akin to the compressible Euler equations with a certain dissipation. In this short note, we construct self-similar solutions of such system in the focusing case for any mass supercritical exponent. To the best of our knowledge these solutions are new, and may formally arise as potential blow-up profiles of the focusing NLS equation. 
 \end{abstract}

\section{Introduction}

\subsection{Historical Background}

 
Consider a complex-valued function $v(x,t)$ solving the  nonlinear Schr\"odinger equation (NLS):
\begin{equation} \label{eq:NLS2}
i\partial_t v(x,t) + \Delta v(x,t) +\mu (v|v|^{p-1})(x,t) = 0,
\end{equation}
where $p$ is an odd integer and $\mu\in \{-1, 1\}$.
Other than the $L^2$ norm of $v$, the Hamiltonian given by
\begin{equation*}\label{E}
E(v) = \frac{1}{2} \int_{\mathbb R^d} | \nabla v |^2  -\mu\frac{1}{p+1} \int_{\mathbb R^d} |v|^{p+1}
\end{equation*}
is also conserved.  Note that if  $\mu=-1$ the equation is called defocusing while if 
$\mu=1$ the equation is focusing. 

The NLS \eqref{eq:NLS2} has a natural scaling, namely  if $v$ is a solution and we set  $v_\lambda (x, t) = \lambda^{2/(p-1)} v(\lambda x, \lambda^2 t)$ then $v_\lambda $  is also a  solution. The critical exponent $s_c$ for which the $\dot H^{s_c}$ norm remains invariant under scaling is given by $s_c = \frac{d}{2} - \frac{2}{p-1}$.

In the defocusing case ($\mu=-1$) it was expected that with enough regularity of the initial data one would always obtain global in time solutions. But in  \cite{Merle-Raphael-Rodnianski-Szeftel:implosion-nls} Merle--Rapha\"el--Rodnianski--Szeftel  proved finite time blow up  for 
certain  defocusing energy supercritical NLS, i. e. when $s_c>1$. In their proof the blow up solution was obtained by using as building blocks smooth, radial self-similar imploding solutions to a compressible Euler equation they had constructed in \cite{Merle-Raphael-Rodnianski-Szeftel:implosion-i}, (see also \cite{Buckmaster-CaoLabora-GomezSerrano:implosion-compressible} by Buckmaster and the first two authors).  In \cite{CaoLabora-GomezSerrano-Shi-Staffilani:nonradial-implosion-schrodinger-TdRd} the authors of this note  generalized the space of functions for which there is blow-up  for an energy  supercritical focusing NLS equation to include also non radial initial data.

 In the defocusing case, the connection mentioned above between NLS and compressible Euler equations occurs via the Madelung transform. In fact,  via this transformation from the equation

\begin{equation} \label{eq:NLSdefocusing}
i\partial_t v(x,t) + \Delta v(x,t) - (v|v|^{p-1})(x,t) = 0,
\end{equation}
one obtains (modulo exponentially small terms in self-similar coordinates) a set of equations equivalent to the compressible Euler ones, namely, for certain values of $\gamma>0$:

\begin{align} \begin{split} \label{eq:CNS}
\rho \p_t u + \rho u \cdot \nabla u &= - \frac{1}{\gamma} \nabla (\rho^{\gamma}) + \nu \Delta u   \\
\p_t \rho + \rm{ div } (\rho u) &= 0.
\end{split} \end{align}

Motivated by the successful work on defocusing NLS mentioned above following this connection between NLS and fluid equations, in this note we use a similar approach to investigate blow-up questions for the focusing NLS
\begin{equation} \label{eq:NLSfocusing}
i\partial_t v(x,t) + \Delta v(x,t) + (v|v|^{p-1})(x,t) = 0.
\end{equation}
The  ultimate goal would be  finding  a highly unstable blow-up profile for the focusing NLS in mass super-critical regime $\frac{d}{2}-\frac{2}{p-1}>0$.

Blow-up questions for focusing NLS equations  have played a fundamental role in the study of dispersive equations, in particular in the last two decades. Here we are not going to give a full survey of the many results that have been proved in terms of blow-up, but we are going to mention only those necessary to frame the result that we are presenting in this note. From a virial argument \cite{Glassey:blowup-nls} one can show that 
blow up can occur when the energy $E$ in \eqref{E} with $\mu=1$ is negative.
There are two types of blow up for NLS,  Type I blow up, which is self-similar, and the Type II which is not. 

For the energy supercritical ($s_c>1$) focusing NLS equation there are no proved self-similar blow up results, the problem remains completely open. For Type II  blow up, in the energy supercritical case,   we recall  \cite{Merle-Raphael-Rodnianski:type-ii-blowup-supercritical-nls}.  If we assume that the problem is mass supercritical but energy subcritical, namely $0<s_c<1$ the picture is more complete. 

Type II blow-up solutions of \eqref{eq:NLSfocusing} were constructed in the mass critical case $s_c=0$ by Merle--Rapha\"el and Perelman in \cite{Merle-Raphael:sharp-upper-bound-blowup-critical-nls,Merle-Raphael:universality-blowup-l2-nls,Merle-Raphael:profiles-quantization-blowup-critical-nls,Merle-Raphael:blowup-upper-bound-critical-nls,Perelman:singularities-critical-nls}. In all these results, all the blow-up profiles constructed are connected to the solitary wave $u(t,x)=e^{iwt}Q_{w}(x)$, where $Q_{w}$ is the unique radial and non negative solution to  the elliptic equation:
\[
\Delta Q_{w}-wQ_{w}+Q_{w}^{p}=0.
\]
In the mass supercritical and energy subcritical case, $0<s_c<1$, the construction of blow-up solutions was carried out through solutions that blow up on a ring or sphere; see, for example, \cite{Fibich-Gavish-Wang:singular-ring-solutions-critical-supercritical-nls}, \cite{Merle-Raphael-Szeftel:collapsing-ring-blowup-mass-supercritical-nls}, \cite{Raphael:existence-stability-blowup-sphere-L2-supercritical-nls},
and \cite{Holmer-Perelman-Roudenko:focusing-d-nls-blowup-sphere}. The Type I blow up in the range $0<s_c<1$ was conjectured in \cite{Sulem-Sulem:focusing-nls-wave-packet},  and \cite{Sulem-Sulem:nls-book} and numerically supported, see for example \cite{Budd:asymptotics-blowup-selfsimilar-nls} and references within. This conjecture has been proved in the range $0<s_c\ll1$  in \cite{Merle-Raphael-Szeftel:stable-selfsimilar-L2-supercritical-nls} ($1\leq d \leq 5$) and  \cite{Bahri-Martel-Raphael:selfsimilar-blowup-slightly-supercritical-nls} ($d\geq 1$), and via a  computer assisted proof in the case $d=p=3$ in \cite{Donninger-Schorkhuber:self-similar-blowup-cubic-nls} and \cite{Dahne-Figueras:self-similar-nls-cgl}. In view of the lack of a complete picture of blow up in this mass supercritical regime it becomes relevant to look for blow up profiles using different methods, such as in this case the connection with  fluid equations.

\subsection{Setup of the problem}


We now present the equation in hydrodynamical variables. Let us perform the Madelung transform and write $v = \sqrt \rho e^{i\psi}$. From \eqref{eq:NLSfocusing}, we have that
\begin{align} \label{eq:NLS_polar}
i \frac{\p_t \rho}{2\sqrt \rho} e^{i\psi} - \p_t \psi \sqrt \rho e^{i\psi} &=
- \left( \sqrt \rho e^{i\psi} \right) \rho^{\frac{p-1}{2}}  - \Delta(\sqrt \rho) e^{i\psi} - \frac{\nabla \rho}{\sqrt \rho} i \nabla \psi e^{i \psi} - \sqrt \rho \div \left( i \nabla \psi e^{i\psi} \right) \\ 
&= 
\sqrt \rho  e^{i \psi} \left(
 - \rho^{\frac{p-1}{2}}  - \frac{\Delta(\sqrt \rho)}{\sqrt \rho} - \frac{\nabla \rho}{\rho} i \nabla \psi  - i  \Delta \psi   + |\nabla \psi |^2 
\right).
\end{align}

Now
\begin{equation*}
\Delta \sqrt \rho = \div \left( \frac{\nabla \rho}{2 \sqrt \rho }\right) = \frac{\Delta \rho}{2 \sqrt \rho} - \frac{1}{4} \frac{|\nabla \rho|^2}{\rho^{3/2}}.
\end{equation*}

Since both $\rho$ and $\psi$ are real, we can identify real and imaginary parts in \eqref{eq:NLS_polar} and obtain
\begin{align} \begin{split} \label{eq:NLS_polar2}
\p_t \psi &= 
\rho^{\frac{p-1}{2}}  + \frac{1}{2\rho} \Delta \rho - \frac{1}{4} \frac{|\nabla \rho|^2}{\rho^2} - |\nabla \psi |^2  \\
\p_t \rho &= 2 \left(
 - \nabla \rho  \cdot \nabla \psi  -   \rho \Delta \psi  \right)
\end{split} \end{align}

In order to see the similarity with the compressible Euler equations, we define $u = \nabla \psi$ and observe that the system takes the form
\begin{align*}
\p_t u &=  \frac{p-1}{p+1} \frac{ \nabla \left( \rho^{\frac{p+1}{2}} \right)}{\rho } - 2 u \nabla u  + \nabla \left( \frac{ \Delta \rho }{2\rho}- \frac{1}{4} \frac{|\nabla \rho|^2}{\rho^2} \right)\\
\p_t \rho &= -2 u \nabla \rho - 2 \rho \div (u) = -2 \div (u \rho).
\end{align*}

We rescale time via $\p_{t'} = 2\p_t$. We have
\begin{align*}
 \p_{t'} u &= \frac{p-1}{2(p+1)} \frac{ \nabla \left( \rho^{\frac{p+1}{2}} \right) }{\rho} - u \nabla u  + \nabla \left( \frac{ \Delta \rho }{4\rho}- \frac{1}{8} \frac{|\nabla \rho|^2}{\rho^2} \right)\\
\p_{t'} \rho &=  -\div (u \rho)
\end{align*}
which corresponds to the compressible Euler equations for $\gamma = \frac{p+1}{2}$ and forcing $\nabla \Delta \sqrt \rho$ 
-- the \textit{quantum pressure}. See the survey \cite{Carles-Danchin-Saut:madelung-gross-pitaevskii-korteweg} by Carles--Danchin--Saut for several aspects of the hydrodynamic formulation of NLS via the Madelung transform.

\subsubsection{Self-similar equations}
We let $\alpha = \frac{\gamma - 1}{2} = \frac{p-1}{4}$. Setting the self-similar ansatz
\begin{align*}
\psi = \frac{ (T-t)^{\frac{2}{r} - 1} }{r} \Psi \left( \frac{x}{(T-t)^{1/r}}, - \frac{\log (T-t)}{r} \right), \\
\rho = \frac{ (T-t)^{\frac{1}{\alpha r} - \frac{1}{\alpha} } }{r} P \left( \frac{x}{(T-t)^{1/r}}, - \frac{\log (T-t)}{r} \right)
\end{align*}
and the variables
\begin{equation*}
s = -\frac{\log(T-t)}{r}, \qquad \mbox{ and } \qquad y = \frac{x}{(T-t)^{1/r}} = x e^s
\end{equation*}
in \eqref{eq:NLS_polar2}, we obtain
\begin{align*}
 -(2-r)\Psi + y \nabla \Psi + \p_s \Psi   &=
r^{-2\alpha + 2}  P^{2\alpha}
- |\nabla \Psi|^2 + (T-t)^{2-\frac{4}{r}} \left( \frac{\Delta P}{2 P} -   \frac{|\nabla P|^2}{4P^2} \right) \\ 
 -\frac{1-r}{\alpha} P + y \nabla P + \p_s P  & =
-2  \nabla P \cdot \nabla \Psi
-2  P \Delta \Psi. \\ 
\end{align*}

Therefore, we obtain
\begin{align} \begin{split} \label{eq:SSequation}
\p_s \Psi   &= (2-r) \Psi -y \nabla \Psi - |\nabla \Psi|^2 + r^{-2\alpha + 2}  P^{2\alpha}
 + e^{(4-2r)s} \left( \frac{\Delta P}{2 P} -   \frac{|\nabla P|^2}{4P^2} \right) \\ 
\p_s P  & = \frac{1-r}{\alpha} P - y\nabla P
-2  \nabla P \cdot \nabla \Psi
-2  P \Delta \Psi. \\ 
\end{split} \end{align}

Defining $S =  2\frac{r^{1-\alpha}}{\sqrt \alpha} P^\alpha$, we have
\begin{align} \begin{split} \label{eq:fulleq_SPsi}
\p_s \Psi &= -(r-2)\Psi - y \nabla \Psi - | \nabla \Psi |^2 + \frac{\alpha}{4} S^2  + e^{(4-2r)s} \frac{\Delta (S^{1/(2\alpha)})}{S^{1/(2\alpha)}} \\ 
\p_s S &= -(r-1)S - y \nabla S -  2\nabla S \cdot \nabla \Psi - 2 \alpha S \Delta \Psi.
\end{split} \end{align}

One can also write the equation in $U, S$ variables, where $U = 2\nabla \Psi$, yielding
\begin{align} \begin{split} \label{eq:fulleq_US}
\p_s U &= -(r-1)U - y \nabla U -  U \cdot \nabla U +  \alpha S \nabla S  + 2e^{(4-2r)s} \frac{\Delta (S^{1/(2\alpha)})}{S^{1/(2\alpha)}} \\ 
\p_s S &= -(r-1)S - y \nabla S -  \nabla S \cdot U -  \alpha S \text{div} (U).
\end{split} \end{align}

We will now look for $U, S$ that are radially symmetric, in a way that we can write $U = U_R \hat{r}$ where $\hat{r}$ is the unitary vector in the radial coordinate, and such that both functions are independent of $s$. Neglecting the exponentially decaying term, after substituting in \eqref{eq:fulleq_US} we obtain the following equations:
\begin{align} \label{eq:SSeq_zeta} \begin{split} 
0 &= -(r-1)U_R - \zeta \p_\zeta U_R - U_R \p_\zeta U_R +  \alpha S \p_\zeta S  \\ 
0 &= -(r-1)S - \zeta \p_\zeta S -  \p_\zeta S  U_R -  \alpha S \p_\zeta U_R - \alpha  (d-1) \frac{SU}{\zeta} .
\end{split} \end{align}


Now we apply the following extra transformation: we define $U_R = \zeta \bar U$ and $S = \zeta \bar S$. We then obtain that 
\begin{align} \begin{split} 
0 &= -(r-1) \bar U - \zeta \p_\zeta \bar U  -  \bar U - \zeta \bar U \p_\zeta  \bar U -  \bar U^2  +  \alpha \bar S^2 + \alpha \zeta \bar S \p_\zeta \bar S  \\ 
0 &= -(r-1) \bar S - \zeta \p_\zeta \bar S - \bar S - \zeta \bar U \p_\zeta \bar S -  \bar U  \bar S -  \alpha \zeta \bar S \p_\zeta \bar U - \alpha   \bar U  \bar S - \alpha (d-1)  \bar U \bar S.
\end{split} \end{align}
Letting $\p_\xi = \zeta \p_\zeta$
\begin{align}  \label{1}
\begin{split} 
\p_\xi \bar U + \bar U \p_\xi \bar U - \alpha \bar S \p_\xi \bar S &= -r \bar U   -  \bar U^2  +  \alpha \bar S^2  \\ 
\p_\xi \bar S + \bar U \p_\xi \bar S + \alpha \bar S \p_\xi \bar U &= -r \bar S    -  \bar U \bar S  - \alpha d  \bar U \bar S .
\end{split} \end{align}
 \color{black}
We can rewrite this system as
\begin{align}
\begin{pmatrix}
1 + \bar U & - \alpha \bar S\\
\alpha \bar S & 1 + \bar U \\
\end{pmatrix} \cdot 
\begin{pmatrix}
\p_\xi \bar U \\
\p_\xi \bar S
\end{pmatrix}
=
\begin{pmatrix}
 -r \bar U   -  \bar U^2  +  \alpha  \bar S^2  \\
-r \bar S    - \bar U \bar S  - \alpha d  \bar U \bar S
\end{pmatrix}.
 \end{align}
The fundamental difference with the defocusing case is that here the matrix is not real diagonalizable.  
One can however still invert the matrix, and using
\begin{align*}
\begin{pmatrix}
1 + \bar U & - \alpha \bar S\\
\alpha \bar S & 1 + \bar U \\
\end{pmatrix}^{-1}
=
\frac{1}{(1+\bar U)^2 + \alpha^2 \bar S^2}
\begin{pmatrix}
1 + \bar U & \alpha \bar S\\
-\alpha \bar S & 1 + \bar U \\
\end{pmatrix}
\end{align*}
and \eqref{1} yields
\begin{align}
\label{dUdS}
\p_\xi \bar U = \frac{N_U}{D}, \qquad \mbox{ and } \qquad \p_\xi \bar S = \frac{N_S}{D}
\end{align}
for
\begin{align}
D &= \left( 1 + \bar U \right)^2 + \alpha^2\bar S^2, \nonumber \\
N_U &= (1 + \bar U ) \left(  -r \bar U   -  \bar U^2  +  \alpha  \bar S^2  \right) + \alpha \bar S \left( -r \bar S    - \bar U \bar S  - \alpha d  \bar U \bar S \right) \nonumber \\
N_S &= -\alpha \bar S \left(   -r \bar U   -  \bar U^2  +  \alpha  \bar S^2  \right) + (1 + \bar U) \left( -r \bar S    - \bar U \bar S  - \alpha d  \bar U  \bar S \right). \label{NuNsD}
\end{align}
Another fundamental difference with the defocusing case is that here $D \neq 0$.

\subsection{Main result}

We now state the main result of this paper:
\begin{theorem} 
Assume that $\alpha d > r-1$, $\alpha=\frac{p-1}{4}$, with $p$ an odd natural number. Then there exists a smooth solution $( U_R(\zeta ),  S(\zeta ))$ to \eqref{eq:SSeq_zeta} with decay at infinity. Moreover, $U_R$ admits a $C^\infty$ odd extension and $S$ admits a $C^\infty$ even extension with respect to $\zeta$, so that the $\mathbb R^d$ radially symmetric vector and scalar fields generated by $U_R$ and $S$, respectively, are smooth. Under the scaling symmetry, we can take $S(0)=1$ and $U'_{R}(0)=-\frac{r-1}{\alpha d}.$\end{theorem}

\begin{figure}[htbp]
    \centering
    \includegraphics[width=0.6\textwidth]{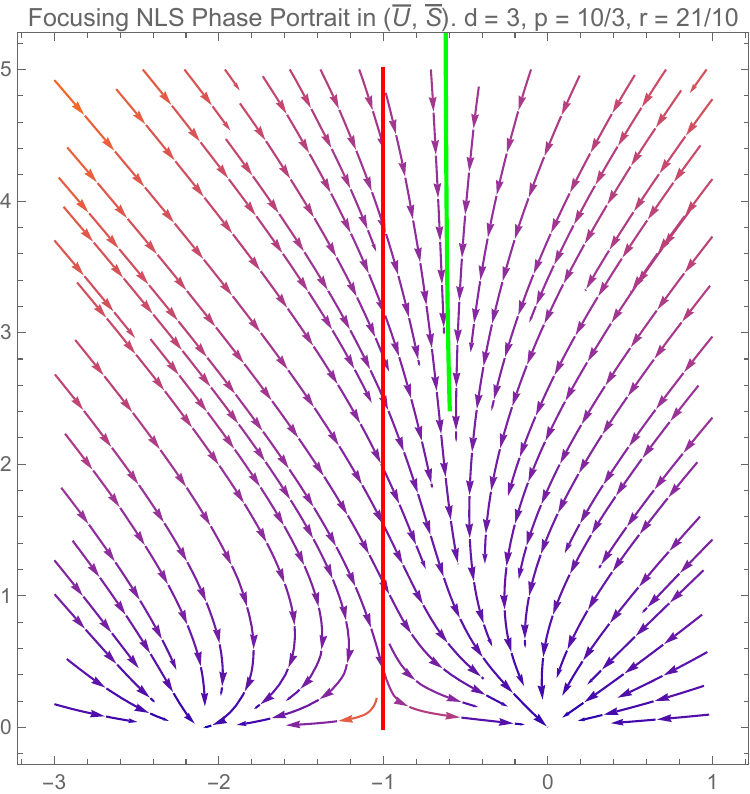}
    \caption{Green: Taylor series expansion starting from $P_0$. The trajectory will reach $(0,0)$. Red: barrier $\bar U = -1$. The solution needs to start from $(\bar U, \bar S) = (\bar U_0, \infty)$ with $\bar U_0 > -1$. There will exist values of $r > 2$ such that this happens if and only if $(p, d)$ lie in the mass supercritical regime.}
      
    \label{fig:nls-focus}
\end{figure}

The proof is postponed to the next section and done in terms of the autonomous system \eqref{dUdS}. It is structured in the following steps:

\begin{enumerate}
    \item We first perform an expansion of the solution as a series around $\xi = -\infty$.
    \item Then, we can bound the coefficients of such expansion to prove its convergence for $\xi < -\log (8)$. At such point standard existence techniques for ODE apply, and we can continue the solution.
    \item In order to prove that the limit of the solution as $\xi \to \infty$ is $(0,0)$, we apply barrier arguments (with the barriers $U = -1$ and $U = 0$, see Figure \ref{fig:nls-focus} with an illustration of the phase portrait) to show that the solution will never cross them and thus the only possibility for such limit is the point $(0,0)$.
\end{enumerate}

\begin{remark}

In order to justify that it is ok to neglect the term $e^{(4-2r)s} \frac{\Delta (S^{1/(2\alpha)})}{S^{1/(2\alpha)}}$, we need to assume that $r > 2$. By choosing $r = 2+\eps$, in terms of $p$, we obtain $(p-1) d > 4 (1+\eps)$, so that for any $p > 1 + 4/d$ we can find such $r$. This corresponds to a \emph{mass supercritical NLS equation}. Indeed as mentioned above one has that
\begin{equation*}
s_c = \frac{d}{2} - \frac{2}{p-1},
\end{equation*}
where the choice of $p_{L^2} = 1 + 4/d$ corresponds $s_c = 0$. Hence our condition $p > p_{L^2}$ corresponds to $s_c > 0$, the  mass supercritical case.

\end{remark}

\begin{remark}
The point $(0, 0)$ is a focus of our phase portrait (cf. Figure \ref{fig:nls-focus}). The rate of decay of any solution approaching it will be dictated by its eigenvalues, concretely if the Jacobian of the focus has eigenvalues $\lambda_1 \leq \lambda_2 < 0$, all solutions will approach at the slower rate $e^{\lambda_2 \xi}$, while two exceptional solutions (forming one smooth curve) will approach at rate $e^{\lambda_1 \xi}$. Such considerations happen to be irrelevant in our case, since for our phase portrait linearized around $(0, 0)$, the Jacobian of the field has a double eigenvalue $\lambda_1 = \lambda_2 = -r$. Thus, any solution (and in particular, our solution), will be $O(e^{-r\xi})$ for large $\xi$. After undoing the change $U_R = e^{\xi} \bar U$ and $S = \bar S e^{\xi}$ the original self-similar solution still decays, in this case at a rate $\frac{1}{R^{(r-1)}}$. This is actually exactly the same rate as the solutions in the defocusing case \cite{Merle-Raphael-Rodnianski-Szeftel:implosion-i}. 
\end{remark}

\begin{remark} 
It is very natural to try to upgrade the radial solution of \eqref{dUdS} into a  solution of \eqref{eq:NLSfocusing} using the techniques developed in \cite{Merle-Raphael-Rodnianski-Szeftel:implosion-i,Buckmaster-CaoLabora-GomezSerrano:implosion-compressible,CaoLabora-GomezSerrano-Shi-Staffilani:nonradial-implosion-compressible-euler-ns-T3R3}. However, the lack of energy estimates and the ill-posedness in Sobolev spaces for the \eqref{eq:fulleq_US} focusing case lead to the need to work in analytic spaces \cite{Metivier:remarks-wellposedness-nonlinear-cauchy}. See \cite{Gerard:remarques-analyse-nls,Thomann:instabilities-supercritical-schrodinger-manifolds} for local existence results in analytic spaces. Another point of view is that with the negative pressure law, the compressible Euler equation is now an elliptic rather than a hyperbolic system \cite[Remark 2.2]{Carles-Danchin-Saut:madelung-gross-pitaevskii-korteweg}. 
\end{remark}

\begin{remark}   
Under the hydrodynamical variables we have that $|v | = \sqrt{\rho}$ and $\nabla v = v\left( \frac{\nabla \rho}{2\rho} + i \nabla \psi\right)$. Therefore, the focusing energy can be reexpressed as
\begin{equation*}
E(\rho, \psi) = \int_{\mathbb R^d} \rho \left( \frac{| \nabla \psi |^2}{2} + \frac{| \nabla \rho |^2}{8 \rho^2 } - \frac{ \rho^{ \frac{p-1}{2} } }{p+1}\right)
\end{equation*}
Taking $t = 0$ in our profiles, the inner parenthesis reads
\begin{align*}
\frac{T^{2/r-2}| \nabla \Psi |^2}{2r^2} + \frac{T^{-2/r}| \nabla P |^2}{8 P^2 } - \frac{ T^{2/r-2} P^{ 2\alpha } }{(p+1)r^{2\alpha}} = T^{2/r-2} \left( \frac{|U|^2}{8r^2} + T^{2-4/r}\frac{| \nabla S |^2}{8 \alpha^2 S^2 } -   \frac{ \alpha S^2 }{4 (p+1)r^{2}}  \right)
\end{align*}

 If we fix the profile to be very close to the solution at the singular time, and using that $r > 2$, the sign of the energy agrees with the sign of
$$\int_0^\infty S^{1/\alpha} \left( \frac{U^2}{2} - \frac{\alpha S^2}{p+1}\right) R^{d-1} dR$$
Numerically, we have observed that this value is negative for $d=3, p=10/3, r=21/10$. 

\end{remark}

\subsection{Acknowledgements}
GCL and GS have been supported by NSF under grant DMS-2052651 and DMS-2306378, and by the Simons Foundation Collaboration
Grant on Wave Turbulence. GCL, JGS and JS have been partially supported by the MICINN (Spain) research grant number PID2021–125021NA–I00.
JGS has been partially supported by NSF under Grants DMS-2245017, DMS-2247537 and DMS-2434314, and by the AGAUR project 2021-SGR-0087 (Catalunya). 
JS has been partially supported by an AMS-Simons Travel Grant. JGS is also thankful for the hospitality of the MIT Department of Mathematics, where parts of this paper were done.

\section{Proof of the Main Theorem}

\subsection{Expansion of the solution as a series}

We consider the following expansion 
\begin{equation}
\bar U = \sum_{ \substack{n = 0 \\ \text{ even } } } U_n e^{n \xi}, \qquad \mbox{ and } \qquad \bar S = \sum_{ \substack{ n=-1 \\ \text{ odd }}} S_n e^{n\xi}. 
\end{equation}
Scaling symmetry allows us to fix $S_{-1} = 1$. In order to deduce the recursion of the coefficients, it is convenient to consider the simplified expressions of $D, N_u, N_s$ that one can obtain simplifying the above expressions (cf. \eqref{NuNsD}):
\begin{align*}
D &= \left( 1 +  \bar U \right)^2 + \alpha^2 \bar S^2, \\
N_U &= -\alpha ^2 d \bar S^2 \bar U-\alpha  (r-1) \bar S^2- \bar U (\bar U+1) (r+\bar U) \\
N_S &= -\alpha ^2 \bar S^3+\bar S \bar U (-\alpha  d+\alpha  r-r-1)+\bar S \bar U^2 (\alpha  (-d)+\alpha -1)-r \bar S.
\end{align*}
Now, note that
\begin{equation*}
D = \sum_{\substack{n = -2 \\ \text{even}}} D_n e^{n \xi}, \qquad 
N_U = \sum_{\substack{n = -2 \\ \text{even}}} N_{U, n} e^{n \xi}, \qquad 
N_S = \sum_{\substack{n = -3 \\ \text{odd}}} N_{S, n} e^{n \xi}.
\end{equation*}
The coefficients are given by
\begin{align*}
D_n &= \delta_{n, 0} + 2 U_n + \sum_{\substack{ i = 0\\ \text{even}}}^n U_i U_{n-i} + \alpha^2 \sum_{\substack{ i = -1 \\ \text{odd} }} S_i S_{n-i} \\ 
N_{U, n} &= - \alpha^2 d \sum_{ \substack{ i = -1 \\ \text{odd }} }^{n+1} \sum_{ \substack{ j = -1 \\ \text{odd }} }^{n-i} S_i S_j U_{n-i-j} - \alpha (r-1) \sum_{ \substack{ i = -1 \\ \text{odd }} }^{n+1} S_i S_{n-i} - \sum_{\substack{ i = 0\\ \text{even }}}^n \sum_{\substack{ j = 0\\ \text{even }}}^{n-i} U_i U_j U_{n-i-j}
- (r+1)\sum_{\substack{ i = 0\\ \text{even }}}^n  U_i U_{n-i} - r U_n\\
N_{S, n} &= - \alpha^2 \sum_{ \substack{ i = -1 \\ \text{odd }} }^{n+2} \sum_{ \substack{ j = -1 \\ \text{odd }} }^{n+1-i} S_i S_j S_{n-i-j} + (-\alpha d + \alpha r -r-1)  \sum_{ \substack{ i = -1 \\ \text{odd }} }^{n} S_i U_{n-i} + (-\alpha d + \alpha -1 ) \sum_{ \substack{ i = -1 \\ \text{odd }} }^{n} \sum_{ \substack{ j = 0 \\ \text{even }} }^{n-i} S_j U_i U_{n-i-j} - rS_n.
\end{align*}

Moreover, we extract the top order terms of the numerators as follows
\begin{equation*}
N_{U, n} = - \alpha^2 d U_{n+2} + \tilde N_{U, n}, \qquad \mbox{ and }\qquad N_{S, n} = -3\alpha^2 S_{n+2} + \tilde N_{S, n}
\end{equation*}
where
\begin{align*}
\tilde N_{U, n} &= - \alpha^2 d \sum_{ \substack{ i = -1 \\ \text{odd }} }^{n+1} \sum_{ \substack{ j = (-1)\vee (-i) \\ \text{odd }} }^{n-i} S_i S_j U_{n-i-j} - \alpha (r-1) \sum_{ \substack{ i = -1 \\ \text{odd }} }^{n+1} S_i S_{n-i} - \sum_{\substack{ i = 0\\ \text{even }}}^n \sum_{\substack{ j = 0\\ \text{even }}}^{n-i} U_i U_j U_{n-i-j}
- (r+1)\sum_{\substack{ i = 0\\ \text{even }}}^n  U_i U_{n-i} - r U_n\\
\tilde N_{S, n} &= - \alpha^2 \sum_{ \substack{ i = -1 \\ \text{odd }} }^{n} \sum_{ \substack{ j = (-1)\vee (-i) \\ \text{odd }} }^{(n+1-i)\wedge n} S_i S_j S_{n-i-j} + (-\alpha d + \alpha r -r-1)  \sum_{ \substack{ i = -1 \\ \text{odd }} }^{n} S_i U_{n-i} + (-\alpha d + \alpha -1 ) \sum_{ \substack{ i = -1 \\ \text{odd }} }^{n} \sum_{ \substack{ j = 0 \\ \text{even }} }^{n-i} S_j U_i U_{n-i-j} - rS_n.
\end{align*}
We have removed the term $S_{-1} S_{-1} U_{n+2} = U_{n+2}$ on the first sum of $N_{U, n}$ and the term $S_{-1}S_{-1} S_{n+2} = S_{n+2}$ (and its permutations) on the first sum of $N_{S, n}$.

Now, we look at the series for $D \p_\xi U$ and $D \p_\xi S$. We have that
\begin{align*}
(D \p_\xi U)_n &= n U_n + 2 \sum_{\substack{i = 0 \\ \text{even }}}^n i U_i U_{n-i} + \sum_{\substack{i = 0 \\ \text{even }}}^n \sum_{\substack{j = 0 \\ \text{even }}}^{n-i} i U_i U_j U_{n-i-j} +\alpha^2 \sum_{\substack{i = 0 \\ \text{even}}}^{n} iU_i \sum_{j = -1}^{n+1-i} S_j S_{n-i-j} + \alpha^2 (n+2)U_{n+2} S_{-1}^2\\
(D \p_\xi S)_n &= n S_n + 2 \sum_{\substack{i = -1 \\ \text{odd }}}^n i S_i U_{n-i} + \sum_{\substack{i = -1 \\ \text{odd }}}^n \sum_{\substack{j = 0 \\ \text{even }}}^{n-i} i S_i U_j U_{n-i-j} + \alpha^2 \sum_{\substack{i = -1 \\ \text{odd}}}^{n}  \sum_{j = (-1)\vee (-i) }^{(n+1-i) \wedge n} iS_i S_j S_{n-i-j} + \alpha^2 n S_{n+2} S_{-1}^2. \\
\end{align*}
Note again that we have isolated the main terms in the right. Defining $\tilde D_{U, n}$ and $\tilde D_{S,n}$ to be the rest and recalling $S_{-1} = 1$, we obtain
\begin{equation*}
(D \p_\xi U)_n = \tilde D_{U, n} + \alpha^2 (n+2) U_{n+2}, \qquad \mbox{ and }\qquad (D \p_\xi S)_n = \tilde D_{S, n} + \alpha^2 n S_{n+2}.
\end{equation*}
Equaling each coefficient of $D\p_\xi U$ with each coefficient of $N_u$ (and the same for $S$) the recurrence can be written as
\begin{equation} \begin{cases} \label{eq:taylor_recursion}
\left( \alpha^2(n+2) + \alpha^2 d \right)U_{n+2} &= \tilde N_{U, n} - \tilde D_{U, n} \\ 
\left(\alpha^2 n + 3\alpha^2 \right)S_{n+2} &= \tilde N_{S, n} - \tilde D_{S, n}
\end{cases}\end{equation}
for every $n \geq -2$. Since $S_{-1} = 1$, and the right hand sides depend at most on coefficients of order $(n+1)$, this allows to solve all the other coefficients as long as the parenthesis of the left hand side are positive since clearly $n+2+ d > 0$ and $n+3>0$.

\subsection{Convergence of the series}
Let us denote $p_n$ to be $1+|U_n|$ for even $n$ and $1+|S_n|$ for odd $n$. Let $n \geq 0$, it is clear that we have
\begin{equation} \label{eq:pre_catalan}
\left| p_{n+2} \right| \leq C \sum_{i, j \geq -1}^{n+2} p_i p_j p_{n-i-j} \mathbbm{1}_{i, j, n-i-j \leq n+1} 
\end{equation}
for some absolute constant $C = C(\alpha, d)$. One can just directly observe this from the expressions of $\tilde N_{U, n}$, $\tilde N_{S, n}$, $\tilde D_{U, n}$, $\tilde D_{S, n}$, since no term with a coefficient $n+2$ appears there. In the cases where the expression is not trilinear, but just bilinear or linear it can be trivially bounded by the trilinear expression above taking $i = 0$ (bilinear) or $i=j=0$ (linear) and noting $p_0 \geq 1$ by definition. Moreover, in the case of $D_{U, n}$, $D_{S, n}$ some factors of $i$ appear, but those are smaller than the factor of $n$ that goes dividing from the companion factor of $U_{n+2}$ or $S_{n+2}$. 

Now, we define $\bar p_n = \max_{i \leq n} p_n$. We take \eqref{eq:pre_catalan} and observe that it cannot happen that two of the three subindex are $-1$ (since then the other would be $n+2$). Therefore, there is at most one of them that is $-1$, and in that case we bound $p_{-1} \leq \bar p_0$. If no subindex is $-1$ we have $i \leq n$ and simply bound $p_i \leq \bar p_{i+1}$. We obtain
\begin{equation*}
\left| \bar p_{n+2} \right| \leq 3C \sum_{\substack{ i, j, k \geq 0 \\ i+j+k = n+1}} \bar p_i \bar p_j \bar p_k 
\end{equation*}
where the factor of $3$ is because one triplet can come from at most three of the previous triplets (one of the three indices is a unit less in the previous sum). 

Therefore, we see that $\bar p_n$ is bounded by the tri-Catalan numbers (also known as number of ternary trees). That is, the sequence defined by $\mathfrak C_0 = 1$ and 
\begin{equation*}
\mathfrak C_{n+1} = \sum_{\substack{ i, j, k \geq 0 \\ i+j+k = n}} \mathfrak C_i \mathfrak C_j \mathfrak C_k.
\end{equation*}
They are known to obey an exponential bound $\mathfrak C_n \les 8^n$. This shows the Taylor series is convergent for $e^{\xi} < 1/8$, After that, we can just use standard ODE existence of solutions to continue our solution $\bar U, \bar S$.

\subsection{Existence of self-similar solutions}

From the convergence of the Taylor series of the previous subsection, we can guarantee that a smooth solution self-similar solution exists for $\xi \in (-\infty, - \log (8) )$ (corresponding to a radial solution on a ball of radius $1/8$  around the origin). Moreover, after undoing the change of variables $U_R = \zeta \bar U = e^\xi \bar U$, and $S = \zeta \bar S = e^{\xi} \bar S$, we obtain that the self-similar solution in $(U, S)$ coordinates is bounded, and that $U$ is odd, while $S$ is even. 

Now, let us bound the global behaviour of the solution. We consider the barrier $\bar U = -1$. The sign of the normal component of the field $(N_u/D, N_s/D)$ to the barrier is just the sign of $N_u$ (given that the normal vector is $(1, 0)$ and $D > 0$). We have that
\begin{equation*}
N_U (-1, S) = \alpha (1 + \alpha d - r) S^2.
\end{equation*}
Therefore, when $\alpha d > r-1$ we have that this sign is positive and therefore solutions cannot traverse $U=-1$ from right to left. 

Similarly, the formula for $\bar U_0$ can be computed from the Taylor recurrence and yields
\begin{equation*}
\bar U_0 = -1 + \frac{\alpha d - (r-1)}{\alpha d}
\end{equation*}
so that $\bar U_0 > -1$ also when $\alpha d > r-1$

Therefore, we conclude that the condition $\alpha d > r-1$ ensures that the solution starts in the region $\bar U > -1$ and moreover cannot exit that region. This condition is particularly interesting since it avoids the possibility of the solution reaching the attractor $(-r, 0)$. On the other hand, this is a if and only if condition, since $\alpha d < r-1$ would imply that the solution starts at $\bar U < -1$ and remains there for all times (making it impossible for the solution to reach $(0, 0)$). In the case $\alpha d = r-1$ one gets the trivial solution $\bar U = -1$, $\bar S = e^{-\xi}$. 

Additionally, we consider the barrier $\bar U = 0$. Here, we have $N_u (0, S) = -S (r + \alpha^2 S^2) < 0$. Therefore, solutions from $\bar U < 0$ cannot cross that barrier. This is the case for our solution, since $\bar U_0 = -(r-1)/(\alpha d) < 0$. We have showed that the condition $\alpha d > r-1$ is required in order for our solution to decay (reach $(0, 0)$), and moreover in such case the solution is bounded on $\bar U \in (-1, 0)$.\\

\bibliographystyle{plain}
\bibliography{references}
\begin{tabular}{l}
  \textbf{Gonzalo Cao-Labora} \\
  {Courant Institute} \\
  {New York University} \\
  {251 Mercer Street, 619} \\
  {New York, NY 10012, USA} \\
  {Email: gc2703@nyu.edu} \\ \\
  \textbf{Javier G\'omez-Serrano}\\
  {Department of Mathematics} \\
  {Brown University} \\
  {314 Kassar House, 151 Thayer St.} \\
  {Providence, RI 02912, USA} \\
  {Email: javier\_gomez\_serrano@brown.edu} \\ \\
  \textbf{Jia Shi} \\
  {Departament of Mathematics} \\
  {Massachusetts Institute of Technology} \\
  {182 Memorial Drive, 2-157} \\
  {Cambridge, MA 02139, USA} \\
  {Email: jiashi@mit.edu} \\ \\
  \textbf{Gigliola Staffilani} \\
  {Departament of Mathematics} \\
  {Massachusetts Institute of Technology} \\
  {182 Memorial Drive, 2-251} \\
  {Cambridge, MA 02139, USA} \\
  {Email: gigliola@math.mit.edu} \\
\end{tabular}
\end{document}